\documentclass[12pt]{amsart}
\usepackage{amsmath,amssymb, euscript, graphicx}
\usepackage[utf8]{inputenc}

\theoremstyle{definition}

\begin{document}

  \title{Preimage Cardinalities of Continuous Functions}
\author {Seljon Akhmedli}

\maketitle

\address{Sheyenne High School, 11th Grade, West Fargo, ND, USA}

\vspace{1cm}

   \begin{center} June 12th, 2020   \end{center}
   
   \bigskip

\begin{center} {\Large ABSTRACT:} We find all subsets of $\mathbb{N} $  which occur as the set of possible cardinalities of preimages of a continuous function. We also study and answer this question for various subclasses of continuous functions. \end{center}

   \section{Introduction}

  It is a well-known fact that there exists no continuous function $f:[0,1]\to \mathbb{R} $ such that every point $x\in \mathbb{R}$ has either two or zero preimages. A somewhat less well-known fact is that there exist continuous functions $f:[0,1]\to \mathbb{R} $ where the preimage $f^{-1}(x)$ has an even cardinality for all $x\in \mathbb{R}$, in fact, one can construct such a function where the set of cardinalities of preimages is $\{0,2,4\}$. 
  
  \medskip
  
 \hspace{4cm}  \includegraphics[width=0.5\textwidth]{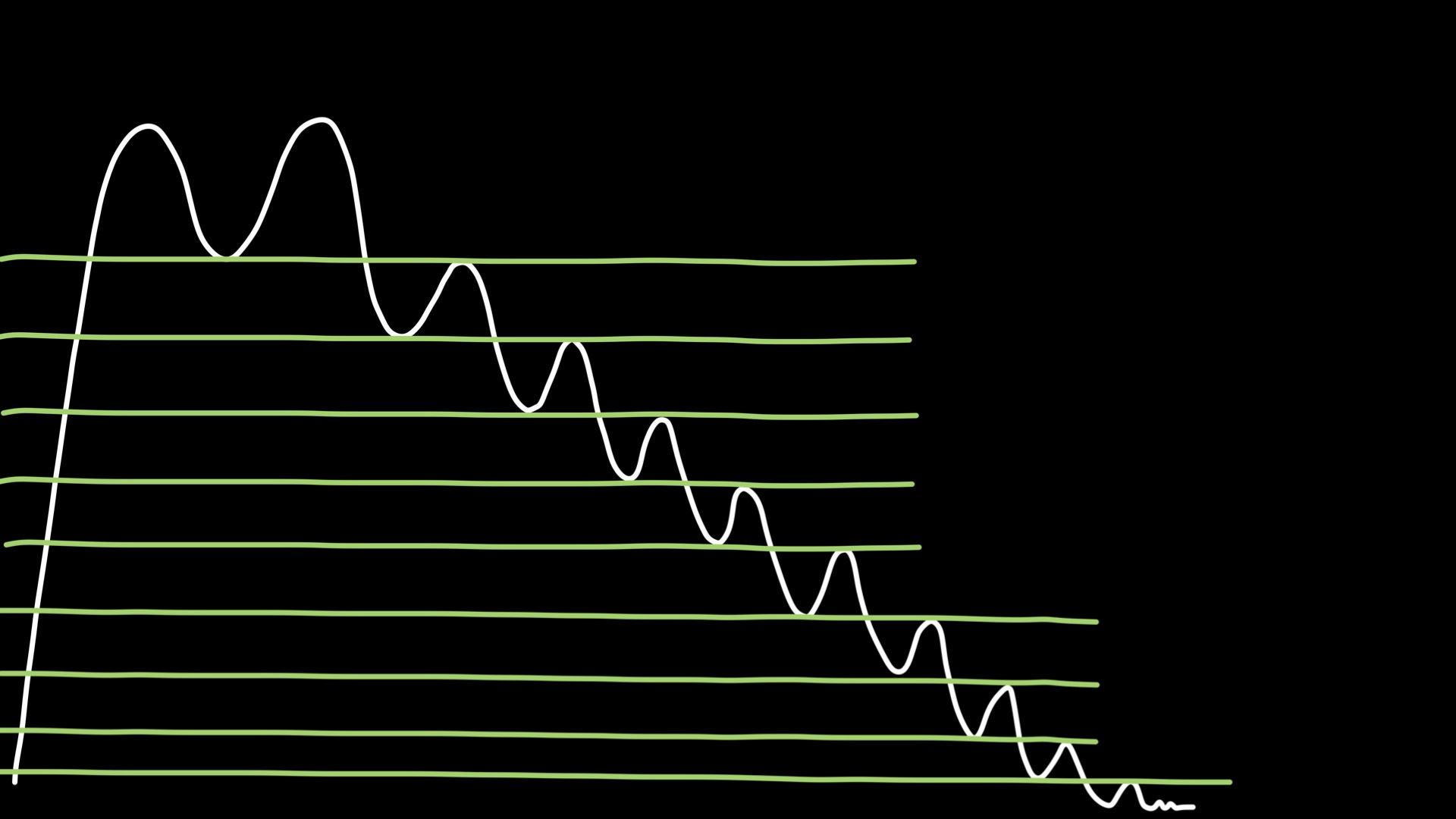}
  
  In this paper, we are interested in the question of what other sets may occur as the set of preimage cardinalities for a continuous function. Given a continuous function $f:[a,b]\to \mathbb{R}, a < b$ we let $$\Omega _f = \{n\in \mathbb{N}\cup \{\infty \} \ : \ \exists x\in \mathbb{R} \ s.t. \ |f^{-1}(x)| = n\}.$$ Our question is what subsets of $\mathbb{N}$ may occur as $\Omega _f$ where $f:[a,b]\to \mathbb{R}$ is a continuous function. Notice that by boundedness of a continuous function in the compact interval, it is necessary that $0\in \Omega _f$. We also study the case when $f$ belongs to a sub-class of the class of continuous functions. 
  
   \medskip
   
  One can also ask what subsets $S$ of $\mathbb{N}\cup \{\infty \} $ occur as $\Omega _f $ but in this question if we let $\infty \in S$ (otherwise, we fall back to the previous question) then it becomes easy: for any non-empty subset $A\subseteq \mathbb{N}$ where $0\in A$, the set $S = A\cup \{\infty \}$ does occur as $\Omega _f$ for some continuous $f:[a,b]\to \mathbb{R}$. 
  
  \medskip
 
  By re-scaling the domain, we may just consider functions over the interval $[0,1]$ in this question (instead of an arbitrary interval). For an arbitrary interval $[a,b]$, $C[a,b]$ will denote the space of all continuous functions $f:[a,b]\to \mathbb{R}$. In this space, we define a subset  $$\mathcal{F}([a,b]) := \{f\in C[a,b] \ : \ f^{-1}(x) \ \mathrm{is \ finite \ for \ all} \ x\in \mathbb{R}\}.$$ In the case of the interval $[0,1]$, we will just write $\mathcal{F}$ instead of $\mathcal{F}([0,1])$. For functions $f\in \mathcal{F}$, we can also write $$\Omega _f = \{n\in \mathbb{N} \ : \ \exists x\in \mathbb{R} \ s.t. \ |f^{-1}(x)| = n\}.$$ If $\Omega _f = S$, then we also say that $f$ is an $S$-function. Notice that a function $f\in \mathcal{F}$ has finitely many local maximums and local minimums, moreover, the end points 0 and 1 are either a local maximum or a local minimum.  
  
  \medskip
  
  Given a subset $S\subseteq \mathbb{N}$ with $|S|\geq 2$, we let $$m_1(S)  = \min (S\backslash \{0\}) \  \mathrm{and} \ m_2(S) = \max S.$$ Here, we use a convention that the maximum of an unbounded subset of natural numbers equals $\infty $. Our main theorem is the following
  
  \medskip
  
  {\bf Theorem 1.} Let $S\subseteq \mathbb{N}$ such that $0\in S$ and $|S|\geq 2$. Then there exists  $f\in C[0,1]$ with $\Omega _f = S$ if and only if $m_2 \geq 2m_1 - 1$. 
  
  \medskip
  
  With slight care, one can indeed make all these functions from the class $C^{\infty }$. An important tool in this construction is the use of functions which "wiggles" infinitely many times. Hence this construction would not work for real analytic functions (by our definition, a function $f:[0,1]\to \mathbb{R}$ is real analytic if it can be extended holomorphically to some domain of $\mathbb{C}$ which includes the unit real segment $[0,1]$; in the sequel, instead of "real analytic" we will simply use the term "analytic"). Motivated by this fact, we also study the same question for all analytic functions, in particular, for polynomials. In fact, by far, not every subset of type $\Omega _f$ for a continuous function occurs as $\Omega _p$ where $p$ is a polynomial. First of all, for a polynomial $p(x)$, we necessarily have $\Omega _p$ is finite (same is true for analytic functions $f:[0,1]\to \mathbb{R}$). But even among finite subsets, for example, the set $\{0, 2, 4\}$ already does not occur as $\Omega _p$. By the following theorem, we achieve a complete characterization of the sets $\Omega _p$ as well where $p(x)$ is an analytic function. To state this theorem, we need to introduce some notations. Let $$A_{+} = \{(n, n+k, n+2k)\in \mathbb{N}^3 \ : \ n\geq 0, k > 0\},$$ \ $$A_{-} = \{(n, n-k, n-2k)\in \mathbb{N}^3 \ : \ n\geq 0, 0 < k\leq \frac{n}{2}\},$$ \ 
  $$B_{+} = \{(n, n+k, n+2k-1)\in \mathbb{N}^3 \ : \ n\geq 0, k \geq 0\},$$ \ $$ B_{-} = \{(n, n-k, n-2k-1)\in \mathbb{N}^3 \ : \ n\geq 0, 0 \leq k\leq \frac{n-1}{2}\},$$ \ 
  $$C_{+} = \{(n, n+k, n+2k-2)\in \mathbb{N}^3 \ : \ n\geq 0, k \geq 0\},$$ \ $$C_{-} = \{(n, n-k, n-2k-2)\in \mathbb{N}^3 \ : \ n\geq 0, 0 \leq k \leq \frac{n}{2}-1\},$$ \ $$A = A_{+}\cup A_{-}, B = B_{+}\cup B_{-}, C = C_{+}\cup C_{-}. $$

  \bigskip
  
  {\bf Theorem 2.} Let $S\subseteq \mathbb{N}$ such that $0\in S$. Then there exists an analytic function $f:[0,1]\to \mathbb{R}$ with $\Omega _f = S$ if and only if there exists a finite sequence $(x_0, x_1, x_2, \dots , x_{2m-1}, x_{2m})$ of natural numbers for some $m\geq 1$ such that the following conditions hold:
  
   i) $x_i\in S$ for all $i\in \{0, \dots, 2m\}$, and for all $z\in S$ there exists $i\in \{0, \dots, 2m\}$ such that $z = x_i$, 
  
   ii) $x_0 = 0$ and $x_{2m} = 0$,
   
   iii) for all $0\leq i\leq m-1, (x_{2i}, x_{2i+1}, x_{2i+2})\in A\cup B\cup C$,
   
   iv) either $|\{i \ : \ 0\leq i\leq m-1, (x_{2i}, x_{2i+1}, x_{2i+2})\in B\}| = 2$ or $|\{i \ : \ 0\leq i\leq m-1, (x_{2i}, x_{2i+1}, x_{2i+2})\in C\}| = 1$,
  
   v) if any triple belongs to $B$, then no triple belongs to $C$.
   
   \medskip
   
   {\bf Remark 1.} If $S$ satisfies conditions i)-v), then the analytic function $p$ with $\Omega _p = S$ is constructed by interpolation method, thus in the statement of the theorem, "analytic" can be replaced with "polynomial".

  \section{Proof of Theorem 1}
  
  In this section we will assume that $S$ is a subset of $\mathbb{N}$ such that $0\in S$ and $|S|\geq 2$. The symbols $I, I_1, I_2, J, J_1, J_2$ will always denote a set $\{1,\dots , n\}$ for some $n\geq 1$ or one of the sets $\mathbb{Z}_{+}, \mathbb{Z}$. Such sets will be called {\em index sets}. We will first consider the special case of Theorem 1 when $1\in S$. In this special case Theorem 1 states the following
   
  \medskip
  
  {\bf Proposition 1.} Let $S\subseteq \mathbb{N}$ such that $0, 1\in S$. Then there exists  $f\in C[0,1]$ with $\Omega _f = S$.
  
  \medskip
  
   To prove Proposition 1, we introduce some terminology. 
   
   {\bf Definition 1.} Let $n\geq 1$. A function $f\in C([a,b])$ is called an {\em $n$-wave function} if there exists real numbers $m < M$ and a sequence $a = x_0 < x_1 < \dots x_{n-1} < x_n = b$ such that $f$ is strictly increasing in the interval $[x_i, x_{i+1}]$ for all even $i\in \{0, \dots , n-1\}$ and strictly decreasing in the interval $[x_i, x_{i+1}]$ for all odd $i\in \{0, \dots , n-1\}$, moreover, $f(x_i) = m$ for all even $i\in \{1, \dots , n\}$ and $f(x_i) = M$ for all odd $i\in \{0, \dots , n\}$. The restriction $f|_{[x_{i-1}, x_i]}$ is called the $i$-th leaf of $f$; the interval $[x_{i-1}, x_i]$ is called the support of this leaf of $f$.
   
   \medskip
   
   Along with $n$-wave functions we will also use the notion of an $n$-stair function.
   
   \medskip
   
    {\bf Definition 2.} Let $n\geq 2$ be a positive odd integer. A function $f\in C([a,b])$ is called an {\em $n$-stair function} if there exists a sequence $a = z_0 < z_1 < \dots z_{n-1} < z_n = b$ such that for all $i\in \{0, \dots , n-1\}$ the restriction $f|_{[z_i, z_{i+1}]}$ is a $k_i$-wave function for some positive odd integer $k_i\geq 1$, and for all $0\leq i < j\leq n-1$, $x\in [z_i, z_{i+1}], y\in [z_j, z_{j+1}]$, $f(x) \leq f(y)$. The sequence $(k_1, \dots , k_n)$ will be called {\em the signature of $f$}. The interval $[z_i, z_{i+1}]$ is called the $(i+1)$-st block of $f$. 
   
   \medskip
   
    \hspace{4cm}  \includegraphics[width=0.5\textwidth]{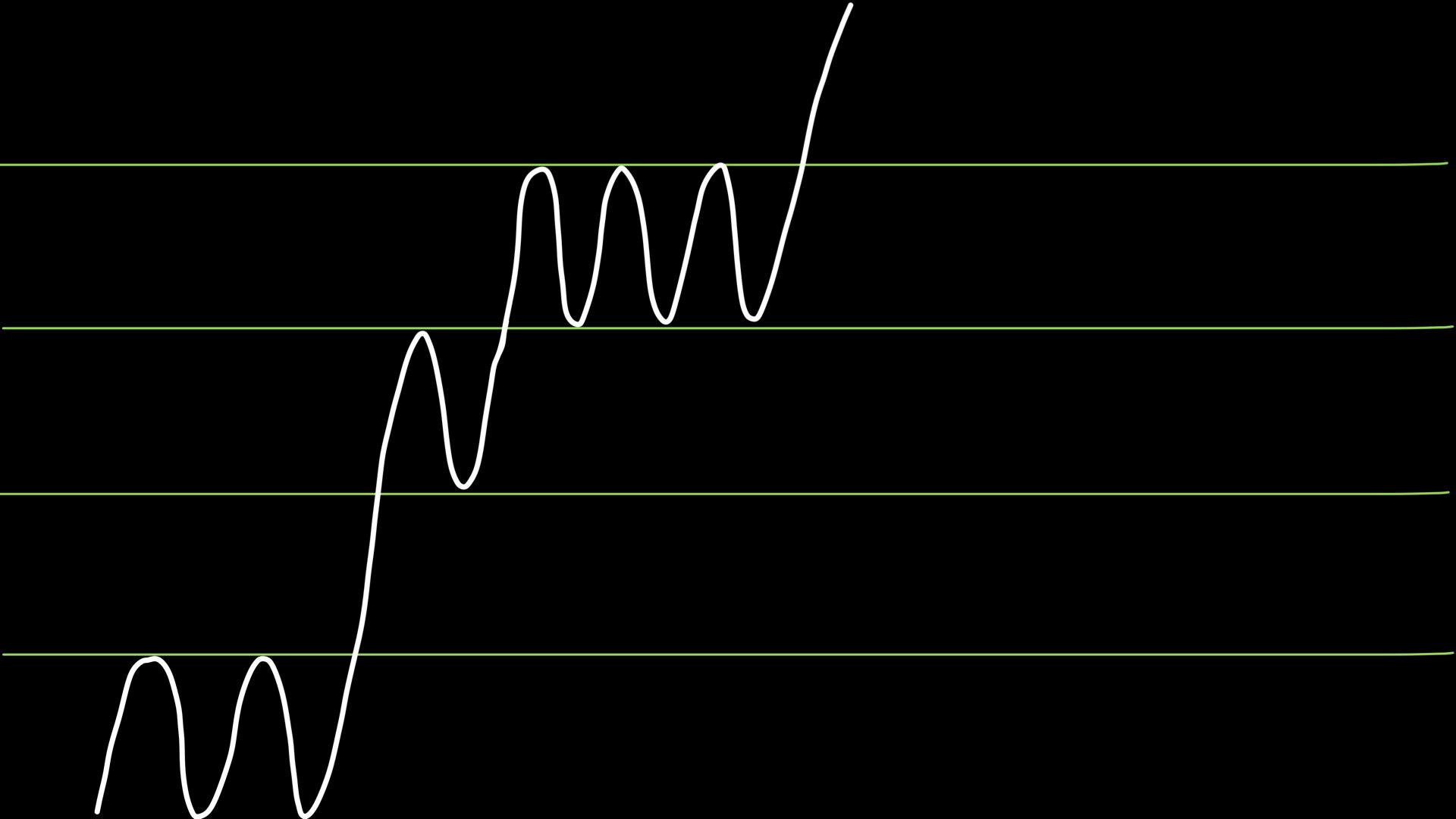}

    In the above picture we have 5-stair function with signature $(5,1,3,7,1)$; so this function has 5 blocks, and the restriction to these blocks are 5-wave, 1-wave, 3-wave, 7-wave and 1-wave functions. 

 \medskip
     
    We also need to define infinite stair functions. This is done similarly as follows
    
    {\bf Definition 3.} A function $f\in C([a,b])$ is called an {\em $\infty $-stair function} if there exists a two-sided sequence $a < \dots < z_{-1} < z_0 < z_1 < \dots  < b$ with $\displaystyle \mathop{\lim }_{n\to \infty }z_{-n} = a$ and $\displaystyle \mathop{\lim }_{n\to \infty }z_n = b$ such that for all $i\in \mathbb{Z}$ the restriction $f|_{[z_i, z_{i+1}]}$ is a $k_i$-wave function for some positive odd integer $k_i\geq 1$, and for all $i < j$, $x\in [z_i, z_{i+1}], y\in [z_j, z_{j+1}]$, $f(x) \leq f(y)$. The sequence $(k_i)_{i\in \mathbb{Z}}$ will be called {\em the signature of $f$}. If $k_i = k$ for all $i\in \mathbb{Z}$ (i.e. if the sequence $(k_i)_{i\in \mathbb{Z}}$ is constant), then we also say that $f$ has a signature $k$. If, for some positive odd $p, q$ and $n_0\in \mathbb{Z}$, $k_i = p$ for all  $i\leq n_0$ and $k_i = q$ for all  $i > n_0$ then we also say that $f$ has a signature $(p,q)$. The intervals $[z_i, z_{i+1}]$ are called {\em the blocks} of $f$. 
    
    \medskip

     \hspace{4cm}  \includegraphics[width=0.5\textwidth]{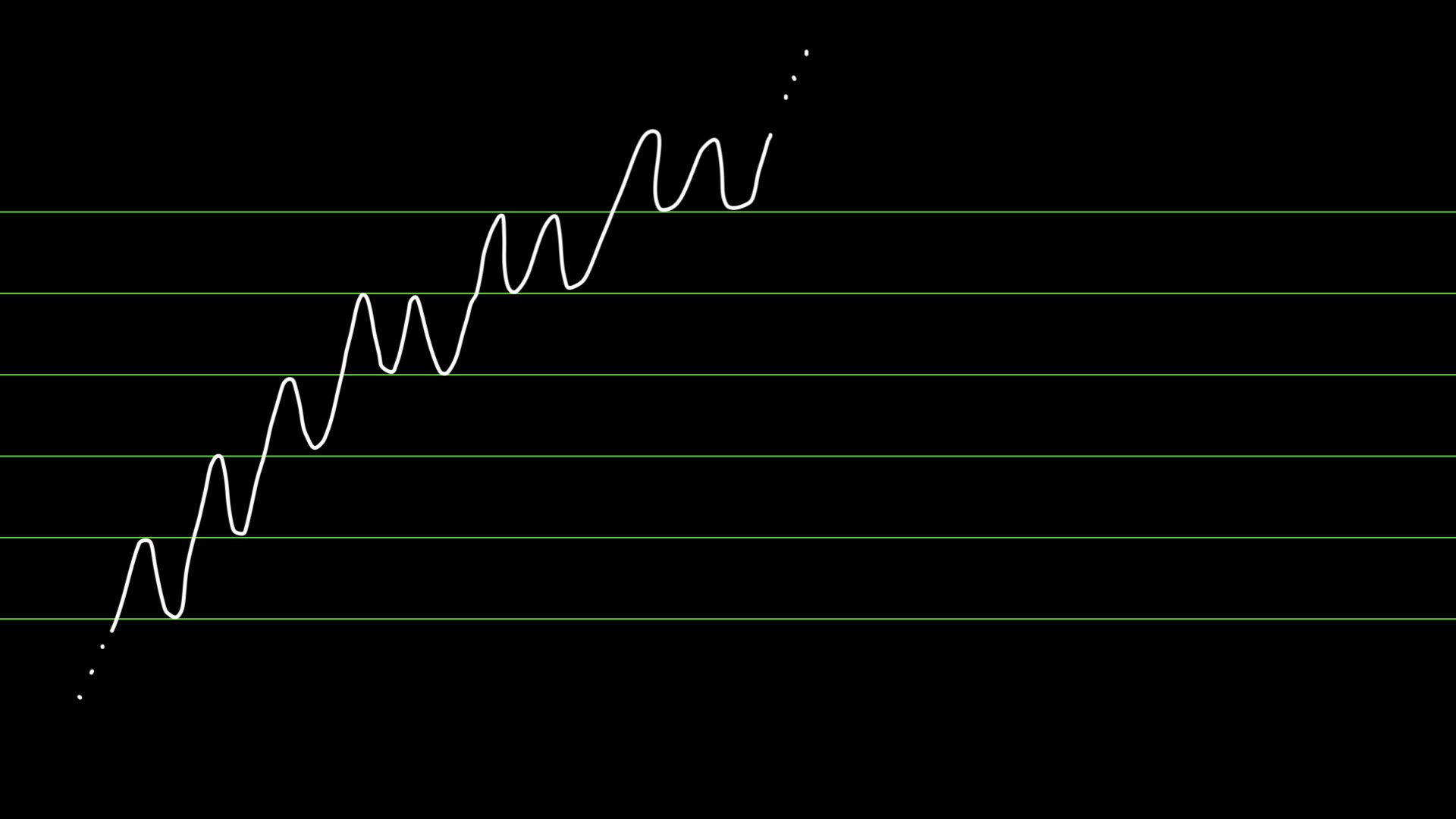}
    
     In the above figure an $\infty $-stair function of signature $(3,5)$ is represented.
     
     \medskip

   {\bf Definition 4.} Let $a\leq c < d\leq b$, $f\in C([a,b]), g\in C([c,d])$ such that $f(c) = g(c)$ and $f(d) = g(d)$. The function $h:[a,b]\to \mathbb{R}$ defined as $$h(x) = \left\{
	\begin{array}{ll}
		f(x)  & \mbox{if } x \in [a,c]\cup [d,b] \\
		g(x) & \mbox{if } x \in [c,d] 
	\end{array}
\right.$$ is called {\em the amendment of $f$ by $g$}. More generally, if $g_i:K_i\to \mathbb{R}, i\in I$ is a collection of continuous functions such that the intervals $K_i = [c_i,d_i]$ have mutually disjoint interiors and $g_i(c_i) = f(c_i), g_i(d_i) = f(d_i)$ for all $i\in I$, then {\em the amendment of $f$ by the collection $(g_i)_{i\in I}$} is defined as $$h(x) = \left\{
	\begin{array}{ll}
		f(x)  & \mbox{if } x \in [a,b]\backslash \displaystyle \mathop{\cup }_{i\in I}K_i \\
		g_i(x) & \mbox{if } x \in K_i, i\in I
	\end{array}
\right.$$
   If $[c,d] = [a,b]$ then the amendment is also called {\em a replacement.}  
   
   \medskip
   
   Notice that the amendment of continuous functions is still continuous. We will be interested in amending $n$-stair functions with $r$-stair functions for various values of $n, r$ (including infinite values) at different sub-intervals.
   
   \bigskip
   
   {\bf Proof of Proposition 1.} In the case $S = \{0,1\}$ the claim is obvious so we may assume that $\{0,1\}$ is a proper subset of $S$. Let $r = \min (S\backslash \{0,1\})$, $(S\backslash \{0,1\})\cap 2\mathbb{N} = \{r_i : i\in I_1\}$ and $(S\backslash \{0,1\})\cap (2\mathbb{N}+1) = \{s_i : i\in I_2\}$. We also define the functions $f:[0,1]\to [0,1], h:[0,1]\to [0,1]$ as $f(x) = x, \forall x\in [0,1]$ and $h(x) = \left\{
	\begin{array}{ll}
		2x  & \mbox{if } x \in [0,\frac{1}{2}] \\
		-x+\frac{3}{2} & \mbox{if } x \in [\frac{1}{2},1] 
	\end{array}
\right.$  We consider the following cases.
   
   \medskip
   
   {\em Case 1.}  $(S\backslash \{0,1\})\cap 2\mathbb{N} = \emptyset $ (i.e. $I_1 = \emptyset $).
   
   Let $C_i, i\in I_2$ be mutually disjoint open intervals in $(0,1)$ such that $\displaystyle \mathop{\cup }_{i\in I_2}\overline{C_i}=[0,1]$. For all $i\in I_2$, we let $g_i:C_i\to C_i$ be a surjective $\infty $-stair function with signature $s_i$. Then we take the identity function $f:[0,1]\to [0,1]$ and for each $i\in I_2$ we amend $f$ with $g_i$ in the interval $C_i$. The amendment $\overline{f}:[0,1]\to \mathbb{R}$ satisfies the conditions $\mathrm{Range}(\overline{f}) = [0,1]$,  $\overline{f}(0) = 0, \overline{f}(1) = 1$ and $\Omega _{\overline{f}} = S$.
   
   \medskip
   
   {\em Case 2.} $2\in S$.
   
 Let $B_i, i\in I_2$ and $C_i, i\in I_1$ be mutually disjoint open intervals in $(0,1)$ such that $\displaystyle \mathop{\cup }_{i\in I_2}\overline{B_i}=[0,\frac{1}{4}]$ and $\displaystyle \mathop{\cup } _{i\in I_1}\overline{C_i}=[\frac{1}{4}, \frac{1}{2}]$.
 
  For all $i\in I_2$ we let $g_i:B_i\to 2B_i$ be a surjective $\infty $-stair function with signature $s_i$ and for all $i\in I_1$ we let $g_i:C_i\to 2C_i$ be a surjective $\infty $-stair function with signature $r_i-1$. Then in the intervals $B_i, i\in I_2$ and $C_i, i\in I_1$ we amend $h$ with $g_i$. The amendment $\overline{h}:[0,1]\to \mathbb{R}$ satisfies the conditions  $\mathrm{Range}(\overline{h}) = [0,1]$, $\overline{h}(0) = 0, \overline{h}(\frac{1}{2}) = 1, \overline{h}(1) = \frac{1}{2}, \overline{h}^{-1}(0) = \{0\}, \overline{h}^{-1}(1) = \{\frac{1}{2}\}, \overline{h}^{-1}(\frac{1}{2}) = \{\frac{1}{4}, 1\}$ and $\Omega _{\overline{h}} = S$.

   \medskip
   
   {\em Case 3.} $r$ is odd, $I_1\neq \emptyset , I_2\neq \emptyset $.
   
   Let $\phi :[0,\frac{1}{2}]\to [0,1]$ be an $\infty $-stair function with signature $(r,r-2)$ such that every block supporting an $r$-wave function lies in $[0,\frac{1}{4}]$ and every block supporting an $(r-2)$-wave function lies in $[\frac{1}{4}, \frac{1}{2}]$, moreover, $|\phi ^{-1}(\frac{1}{2})| = r-1$. Let $B_i, i\in J_2$ and $C_i, i\in J_1$ be mutually disjoint open blocks of $\phi $ such that $\displaystyle \mathop {\cup }_{i\in J_2}\overline{B_i}=[0,\frac{1}{4}]$ and $\displaystyle \mathop{\cup } _{i\in J_1}\overline{C_i}=[\frac{1}{4}, \frac{1}{2}]$.
   
   \medskip
   
   Let $u_i, i\in J_2$ be a leaf of $\phi $ in the block $B_i$, and $v_i, i\in J_1$ be a leaf of $\phi $ in the block $C_i$. Let also $a:J_1\to I_1$ and $b:J_2\to I_2$ be surjective maps (such maps exist because the sets $J_1$ and $J_2$ are both infinite). We replace $u_i, i\in J_2, v_i, i\in J_1$ respectively with $\infty $-stair functions $\overline{u_i}$ of signature $s_{b(i)}-r+1$ and $\overline{v_i}$ of signature $r_{a(i)}-r$. Then we amend $f$ with functions $\overline{u_i}, i\in J_2$ and $\overline{v_i}, i\in J_1$.
   
   \medskip
   
   {\em Case 4.} $r$ is even, $r\geq 4$, $I_1\neq \emptyset , I_2\neq \emptyset $.
   
   This case is very similar to the previous case with the only difference that $\phi :[0,1]\to [0, \frac{1}{2}]$ is assumed to be an $\infty $-stair function with signature $(r+1,r-3)$ such that every block supporting an $(r+1)$-wave function lies in $[0,\frac{1}{4}]$ and every block supporting an $(r-3)$-wave function lies in $[\frac{1}{4}, \frac{1}{2}]$, $|\phi ^{-1}(\frac{1}{2})| = r-1$, and the signatures of $\overline{u_i}, \overline{v_i}$ are $s_{b(i)}-r$ and $r_{a(i)}-r+3$ respectively.
   
   \medskip
   
   {\em Case 5.} $(S\backslash \{0,1\})\cap (2\mathbb{N}+1) = \emptyset $ (i.e. $I_2 = \emptyset $).
   
   \medskip
   
   Since $S\backslash \{0,1\}$ does not contain an odd element we have that $r$ is even. Let $H:[-\frac{1}{2},1]\to [0,1]$ defined as $h(x) = \left\{
	\begin{array}{ll}
		-x  & \mbox{if } x \in [-\frac{1}{2},0] \\
		h(x) & \mbox{if } x \in [0,1] 
	\end{array}
\right.$
   
  and $\psi :[0,\frac{1}{2}]\to [0,1]$ be an $\infty $-stair function with signature $(r-1,r-3)$ such that $|\psi ^{-1}(\frac{1}{2})| = r-2.$ Let also $D_i, i\in J$ be mutually disjoint open blocks of $\psi $ such that $\displaystyle \mathop{\cup }_{i\in J}\overline{D_i}=[\frac{1}{4},\frac{1}{2}]$. Let $w_i, i\in J$ be a leaf of $\psi $ in the block $D_i$, and $c:J\to I_1$ be a surjective function. We replace $w_i, i\in J$ with $\infty $-stair function $\overline{w_i}$ of signature $r_{c(i)}-r+3$. Then we amend $H$ with functions $\overline{w_i}, i\in J$. The amendment $\overline{H}:[-\frac{1}{2},1]\to \mathbb{R}$ satisfies the conditions $\mathrm{Range}(\overline{H}) = [0,1]$, $\overline{H}(-\frac{1}{2}) = \overline{H}(1) = \frac{1}{2}, \overline{H}(0) = 0, \overline{H}(\frac{1}{2}) = 1$ and $\Omega _{\overline{H}} = S$. By re-scaling the domain of $\overline{H}$ we can obtain a function $F:[0,1]\to \mathbb{R}$ with $\Omega _F = S$.  \ $\square $

  \bigskip
  
  For the proof of Theorem 1, we need to introduce the following

  {\bf Definition 5.} Let $S\subseteq \mathbb{N}$ and $0\in S$. $S$ is called of type I if for all even $x\in S$ there exists an odd $y\in S$ such that $y > x$. If $S$ is not of type I then we say that $S$ is of type II. 
  
  \medskip
  
   The above definition implies that if $S$ is finite, then $S$ is of type II if and only if $\max S$ is even.

   \bigskip
   
   Now we are ready for the proof of the main theorem.
   
   \medskip
   
   {\bf Proof of Theorem 1.} First we will prove the "only if part". 
   
   \medskip
   
     If $(0,1)$ does not contain any local maximum point of $f$ then either $f$ is a strictly monotone function or $f$ has only one local minimum in $(0,1)$ and local maximums at 0 and 1. In the latter case $\Omega _f = \{0,1,2\}$ whereas in the former case $\Omega _f = \{0,1\}$; the inequality  $m_2\geq 2m_1-1$ holds in both cases.
    
    \medskip
    
     Thus we may assume that $(0,1)$ does contain a local maximum point. Similarly, we may assume that $(0,1)$ contains a local minimum point. Then, without loss of generality we may also assume that at least one of the ends 0, 1 is not a local maximum. Let $y_1$ be the biggest value of $f$ in $[0,1]$  such that $f^{-1}(y_1)$ contains a local maximum in $(0,1)$. If we have $l$ local maximums in $f^{-1}(y_1)\cap (0,1)$, then $|f^{-1}(y_1)| \leq l+1$, and for a sufficiently small $\epsilon > 0$, $|f^{-1}(y_1-\epsilon )| \geq 2l$; moreover, if $|f^{-1}(y_1)| = l+1$, then $|f^{-1}(y_1-\epsilon )| \geq 2l+1$. Thus $|f^{-1}(y_1-\epsilon )| \geq 2|f^{-1}(y_1)|-1$. Hence we obtain that $m_2\geq 2m_1-1$.

    The "if part" is proved in Proposition 1 when $1\in S$. So we will assume that $1\notin S$.  Let $S\cap 2\mathbb{Z} = \{r_i \ : \ i\in I_1\}$ and $S\cap (2\mathbb{Z}+1) = \{s_i \ : \ i\in I_2\}$. We also let $p = \min (S\backslash \{0\})$. Then $p\geq 2$.

     We will consider the following cases.

    {\em Case 1.} $S$ is of type I.
    
    We can find an odd $q\in S$ such that $q\geq 2p-1$. Let $C_i, i\in \mathbb{Z}$ be mutually disjoint open intervals such that $\displaystyle \mathop{\cup }_{i\in \mathbb{Z}}\overline{C_i} = [\frac{1}{4}, \frac{3}{4}]$ and $C_i$ is on the left of $C_j$ for all $i<j$ (i.e. for all $x\in C_i, y\in C_j, x < y$). 
    
    We let $f:[\frac{1}{4}, \frac{3}{4}]\to [\frac{1}{4}, \frac{3}{4}]$ be a surjective $\infty $-stair function of signature $(\dots , 1, 2p-1, 1, 2p-1, 1, 2p-1, \dots )$ such that every block supported on $C_i$ for even $i\in \mathbb{Z}$ is a 1-wave function, and every block supported on $C_i$ for odd $i\in \mathbb{Z}$ is a $(2p-1)$-wave function. We also let $f_1:[0, \frac{1}{4}]\to [0, \frac{1}{4}]$ and $f_2:[\frac{3}{4},1]\to [\frac{3}{4},1]$ be surjective $(2p-1)$-wave functions. 
    
    Let $u_j$ be leaf of $f_j, j\in \{1,2\}.$ Let also $v_i$ be a leaf of $f$ supported on the block $C_{2i+1}, i\in \mathbb{Z}$ and $w_i = f|_{C_{2i}}, i\in \mathbb{Z}$. So $w_i$ is a 1-wave function for all $i\in \mathbb{Z}$. 
    For all $i\in \mathbb{Z}$, we replace $v_i$ with an $\infty $-stair function $\overline{v_i}$ of signature $(q-2p+2)$; for all $j\in \{1,2\}$ we also replace $u_j$ with an $\infty $-stair function $\overline{u_j}$ of signature $(q-2p+2)$. Furthermore, for all $i\in 2\mathbb{Z}$ we replace $w_i$ with an $\infty $-stair function $\overline{w_i}$ of signature $(r_{a(i)}-1,r_{a(i)}+1)$, and for all $i\in (2\mathbb{Z}+1)$ we replace $w_i$ with an $\infty $-stair function $\overline{w_i}$ of signature $(s_{b(i)})$ where $a:2\mathbb{Z}\to I_1, b:2\mathbb{Z}+1\to I_2$ are surjective maps.  
    The blocks of $\overline{w_i}, i\in 2\mathbb{Z}$ are either $(r_{a(i)}-1)$-wave functions or $(r_{a(i)}+1)$-wave functions. For each $i\in 2\mathbb{Z}$, in the former case we replace a leaf of it with $\infty $-stair function of signature $(q-r_{a(i)}+2)$ whereas in the latter case e replace a leaf of it with $\infty $-stair function of signature $(q-r_{a(i)})$. Let $\overline{\overline{w_i}}$ be this amendment of $\overline{w_i}$. 
    
    \hspace{4cm}  \includegraphics[width=0.5\textwidth]{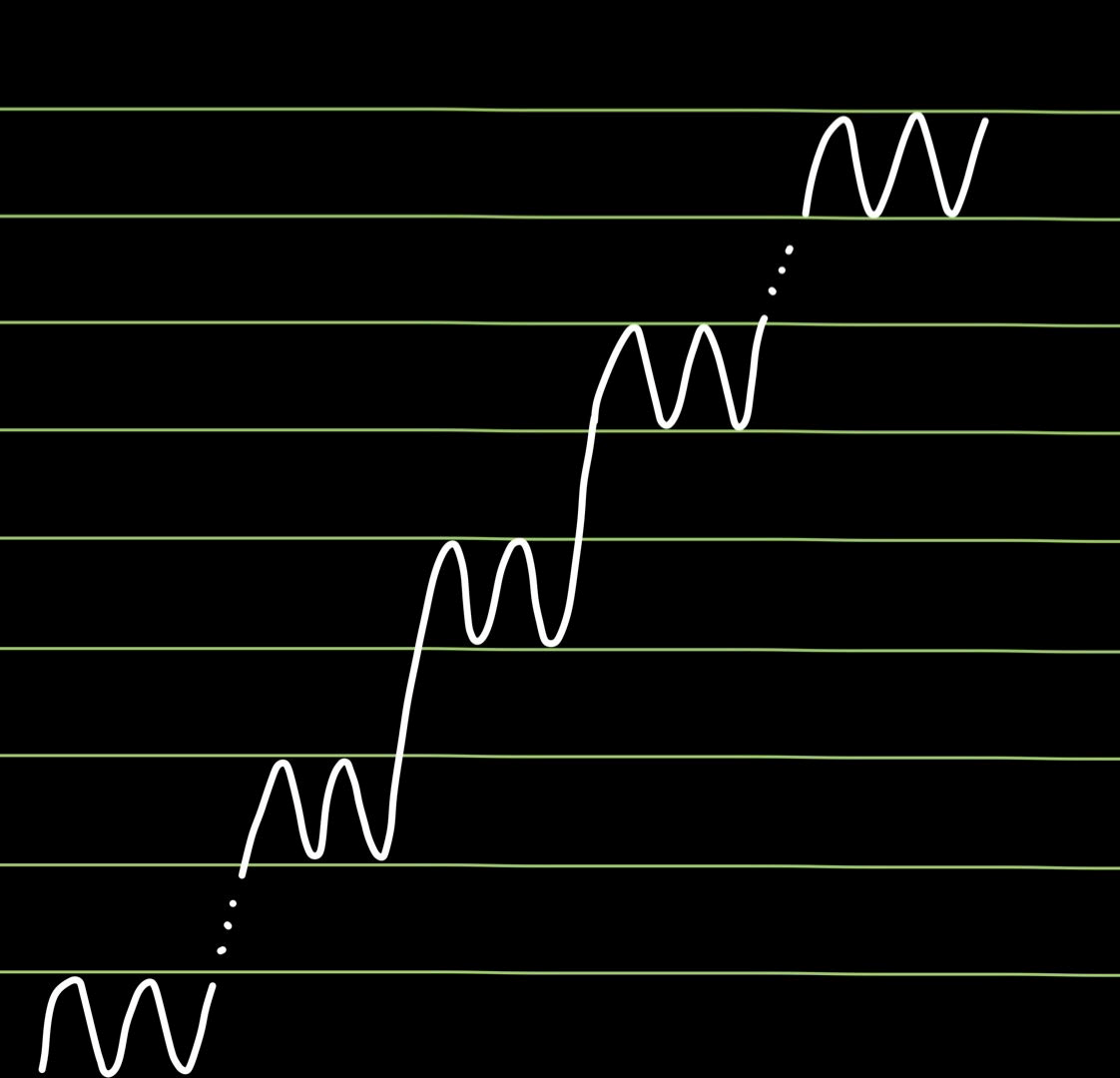}
    
    Now we amend $F:[0,1]\to [0,1]$ with functions $\overline{u_j}, j\in \{1,2\}$ and $\overline{v_i}, i\in \mathbb{Z},  \overline{w_i}, i\in 2\mathbb{Z}+1$ and  $\overline {\overline{w_i}}, i\in 2\mathbb{Z}$ where $F$ is defined as $F(x) = \left\{
	\begin{array}{lll}
		f_1(x)  & \mbox{if } x \in [0, \frac{1}{4}] \\
		f(x) & \mbox{if } x \in [\frac{1}{4},\frac{3}{4}] \\
		f_2(x)  & \mbox{if } x \in [\frac{3}{4},1]
	\end{array}
\right.$  (The above picture represents $F$ when $p=3$).

  For the amendment $\overline{F}:[0,1]\to \mathbb{R}$ we have $\Omega _{\overline{F}} = S$.

  \medskip
  
   {\em Case 2:} $S$ is of type II.
   
   Let the intervals $C_i, i\in \mathbb{Z}$ be as in Case 1 and $a:2\mathbb{Z}\to I_1, b:2\mathbb{Z}+1\to I_2$ are surjective maps. We can find an even $q\in S$ such that $q-1\geq 2p-1$. 
   
   We let $h:[\frac{1}{4}, \frac{3}{4}]\to [\frac{1}{4}, \frac{3}{4}]$ be a surjective $\infty $-stair function of signature $(\dots , 1, 2p-1, 1, 2p-1, 1, 2p-1, \dots )$ such that every block supported on $C_i$ for even $i\in \mathbb{Z}$ is a 1-wave function, and every block supported on $C_i$ for odd $i\in \mathbb{Z}$ is a $(2p-1)$-wave function. We also let $h_1:[0, \frac{1}{4}]\to [0, \frac{1}{4}]$ be a surjective $(2p-1)$-wave function and $h_2:[\frac{3}{4},\frac{7}{8}]\to [\frac{3}{4},\frac{7}{8}]$ be a surjective $(2p-2)$-wave function with $h_2(\frac{7}{8}) = \frac{3}{4}$.
   
   Let $u_j$ be the first leaf of $h_j, j\in \{1,2\}$, and $v_i$ be a leaf of $h$ supported on the block $C_{2i+1}, i\in \mathbb{Z}$ and $w_i = h|_{C_{2i}}, i\in \mathbb{Z}$.
   
   Let the last two leafs of $h_2$ be supported on the interval $D = (d_1, d_2)$. We take a function $\xi : D\to [\frac{3}{4},\frac{7}{8}]$ such that $\xi ^{-1}(\frac{3}{4}) = \{d_1,d_2\}, \xi ^{-1}(\frac{7}{8}) = \{d_3, d_4\}$ for some $d_3, d_4\in D$, (so $|\xi ^{-1}(\frac{3}{4})| = |\xi ^{-1}(\frac{7}{8})| = 2$) and $|\xi ^{-1}(y)| = 4$ for all $y\in (\frac{3}{4},\frac{7}{8})$. Then we amend $h_2$ with $\xi $; let $h_3$ be this amendment. 
   
    For all $i\in \mathbb{Z}$, we replace $v_i$ with an $\infty $-stair function $\overline{v_i}$ of signature $(q-2p+1)$. We also amend $u_1$ with an $\infty $-stair function $\overline{u_1}$ of signature $(q-2p+1)$, and $u_2$ with an $\infty $-stair function $\overline{u_2}$ of signature $(q-2p+1)$. Finally, for all $i\in 2\mathbb{Z}$ we replace $w_i$ with an $\infty $-stair function $\overline{w_i}$ of signature $(s_{a(i)}-2,s_{a(i)})$, and for all $i\in (2\mathbb{Z}+1)$ we replace $w_i$ with an $\infty $-stair function $\overline{w_i}$ of signature $(r_{b(i)}-1)$.
    
    The blocks of $\overline{w_i}, i\in 2\mathbb{Z}$ are either $(s_{a(i)}-2)$-wave functions or $(s_{a(i)})$-wave functions. For each $i\in 2\mathbb{Z}$, in the former case we replace a leaf of it with $\infty $-stair function of signature $(q-s_{a(i)}+2)$ whereas in the latter case e replace a leaf of it with $\infty $-stair function of signature $(q-s_{a(i)})$. Let $\overline{\overline{w_i}}$ be this amendment of $\overline{w_i}$. 
    
    Now let $H:[0,1]\to [0,1]$ be defined as $H(x) = \left\{
	\begin{array}{llll}
		h_1(x)  & \mbox{if } x \in [0, \frac{1}{4}] \\
		h(x) & \mbox{if } x \in [\frac{1}{4},\frac{3}{4}] \\
		h_3(x)  & \mbox{if } x \in [\frac{3}{4},\frac{7}{8}] \\
		-6x+6 & \mbox{if } x \in [\frac{7}{8},1]
	\end{array}
\right.$

  We amend $H$  with functions $\overline{u_j}, j\in \{1,2\}$ and $\overline{v_i}, \overline{w_i}, i\in 2\mathbb{Z}+1$ and  $\overline {\overline{w_i}}, i\in 2\mathbb{Z}$. For the amendment $\overline{H}:[0,1]\to \mathbb{R}$ we have $\Omega _{\overline{H}} = S$.

  \ $\square $

    \bigskip
    
    {\bf Remark 2.} In the proof of Proposition 1 and Theorem 1, we use  $\infty $-stair functions $\sigma :[a,b]\to [c,d]$. We can allow "the heights of the waves" to be small enough and fast converging to zero so we can build such a function $\sigma _1:[a,b]\to [c,c+\delta ]$ of class $C^{\infty }$, where $\delta $ is some (possibly small) positive number and $\sigma ^{(n)}(a) = \sigma ^{(n)}(b) = 0$ for all $n\geq 1$. Then we re-scale the range by defining $\sigma _2(x) = \frac{c}{\delta }(\delta + c - d)+\frac{d-c}{\delta }\sigma _1(x)$ thus make $\sigma _2$ of class $C^{\infty }$. This allows to make the function $f$ in the claim of Proposition 1 and Theorem 1 from the $C^{\infty }$ class.

    \bigskip
    
    \section{Proof of Theorem 2}
    
    Let us first recall a well-known fact that if an analytic function has infinitely many zeros on a compact interval then it is identically zero. Thus an analytic function $p:[0,1]\to \mathbb{R}$ has finitely many zeros and  finitely many critical points. It also belongs to $\mathcal{F}$. Then there exists a finite sequence $y_1 < y_2 < \dots < y_{m-1}$ such that any local maximum or local minimum point of $p$ belongs to $\displaystyle \mathop{\sqcup }_{1\leq i\leq m-1}p^{-1}(y_i)$, and each of the subsets $p^{-1}(y_i)$ contains at least one local maximum or local minimum. Let $y_0 = -\infty $ and $y_m = \infty $, and let $\phi :\mathbb{R} \to \mathbb{N}$ be defined as $\phi (y) = |p^{-1}(y)|$. Then in any interval $(y_i, y_{i+1}), 0\leq i\leq m-1$, $\phi $ is constant. Let $z_i$ be a point in $(y_i, y_{i+1}), 0\leq i\leq m-1$.
    
    \medskip

    We let $x_{2i} = |p^{-1}(z_i)|, 0\leq i\leq m$,  and $x_{2i+1} = |p^{-1}(y_i)|, 1\leq i\leq m$. Then for all $0\leq i\leq m-1$, 
    
    if $p^{-1}(y_i)$ does not contain end points, contains $l$ local minimums and $r$ local maximums, and $r < l$, then  $(x_{2i}, x_{2i+1}, x_{2i+2})\in A_{+}$;
    
    if $p^{-1}(y_i)$ does not contain end points, contains $l$ local minimums and $r$ local maximums, and $r > l$, then  $(x_{2i}, x_{2i+1}, x_{2i+2})\in A_{-}$;
     
    if $p^{-1}(y_i)$ contains one end point as a local minimum, $l$ local minimums and $r$ local maximums, and $r \leq l$, then  $(x_{2i}, x_{2i+1}, x_{2i+2})\in B_{+}$;
      
    if $p^{-1}(y_i)$ contains one end point as a local minimum, $l$ local minimums and $r$ local maximums, and $r > l$, then  $(x_{2i}, x_{2i+1}, x_{2i+2})\in B_{-}$;
      
    if $p^{-1}(y_i)$ contains one end point as a local maximum, $l$ local minimums and $r$ local maximums, and $r \leq l$, then  $(x_{2i}, x_{2i+1}, x_{2i+2})\in B_{+}$;
     
    if $p^{-1}(y_i)$ contains one end point as a local maximum, $l$ local minimums and $r$ local maximums, and $r > l$, then  $(x_{2i}, x_{2i+1}, x_{2i+2})\in B_{-}$;
     
      if $p^{-1}(y_i)$ contains two end points of opposite types, $l$ local minimums and $r$ local maximums, and $r \leq l$, then  $(x_{2i}, x_{2i+1}, x_{2i+2})\in C_{+}$;
      
      if $p^{-1}(y_i)$ contains two end points of opposite types, $l$ local minimums and $r$ local maximums, and $r > l$, then  $(x_{2i}, x_{2i+1}, x_{2i+2})\in C_{-}$;
      
      if $p^{-1}(y_i)$ contains two end points both of which are local maximums, $l$ local minimums and $r$ local maximums, and $r > l+2$, then  $(x_{2i}, x_{2i+1}, x_{2i+2})\in C_{-}$;
      
      if $p^{-1}(y_i)$ contains two end points both of which are local maximums, $l$ local minimums and $r$ local maximums, and $r \leq l+2$, then  $(x_{2i}, x_{2i+1}, x_{2i+2})\in C_{+}$;
      
      if $p^{-1}(y_i)$ contains two end points both of which are local minimums, $l$ local minimums and $r$ local maximums, and $l \geq r+2$, then  $(x_{2i}, x_{2i+1}, x_{2i+2})\in C_{+}$;
      
      if $p^{-1}(y_i)$ contains two end points both of which are local minimums, $l$ local minimums and $r$ local maximums, and $l > r+2$, then  $(x_{2i}, x_{2i+1}, x_{2i+2})\in C_{-}$;
      
      in all other cases $x_{2i} = x_{2i+1} = x_{2i+2}$. 
      
      \medskip
      
      Thus we proved the "only if part" of the theorem. To prove the "if part", we first observe the following simple fact.
      
      \medskip
      
      {\bf Lemma 1.} Let $f\in \mathcal{F}([a,b])$ and $y_1 < y_2 < y_3$ such that for all $y \in [y_1, y_3]\backslash \{y_2\}$, the set $f^{-1}(y)$ does not contain any local maximum or local minimum, and $|f^{-1}(y_1)| = n, |f^{-1}(y_2)| = n+k, |f^{-1}(y_3)| = n+2k$ for some $n\geq 0, k\in \mathbb{Z}\backslash \{0\}$. Let also $r, s\in \mathbb{Z}\backslash \{0\}$ such that $r+s = k$. Then $f$ can be amended in the set $f^{-1}([y_1, y_3])$ such that the amendment $F$ belongs to $\mathcal{F}([a,b])$, moreover, there exists $z_1, z_2 \in (y_1, y_3)$ such that $z_1 < y_2 < z_2$, for all $y \in [y_1, y_3]\backslash \{z_1, z_2, y_2\}$, the set $F^{-1}(y)$ does not contain any local maximum or local minimum, and $|F^{-1}(y_1)| = n, |F^{-1}(z_1)| = n+r, |F^{-1}(y_2)| = n+2r, |F^{-1}(z_2)| = n+2r+s, |F^{-1}(y_3)| = n+2k$.  $\square $

      \medskip
      
      Now, given any sequence $(x_0, x_1, x_2, \dots , x_{2m-1}, x_{2m})$ satisfying conditions i)-v) we have two cases: 1) None of the triples $(x_{2i}, x_{2i+1}, x_{2i+2}), 0\leq i\leq m-1$ belongs to $C$ or 2) Exactly one of the triples $(x_{2i}, x_{2i+1}, x_{2i+2}), 0\leq i\leq m-1$ belongs to $C$.
      
      In the former case we have $0\leq p < q\leq m-1$ such that $$(x_{2p}, x_{2p+1}, x_{2p+2}), (x_{2q}, x_{2q+1}, x_{2q+2})\in B$$ and in the latter case we have $0\leq k\leq m-1$ such that $(x_{2k}, x_{2k+1}, x_{2k+2})\in C$.

      We can take a sequence $(0, 1, \dots , m-1)$ and inductively build a piecewise linear function $f:[0,1]\to \mathbb{R}$ (i.e. $f$ is continuous, has finitely many singularity points, and in between these singularities it is defined as a polynomial of degree at most 1, such that 
      
      (i) if $c$ is a local maximum or local minimum, then $f(c)\in \{0, 1, \dots , m-1\}$
      
      (ii) $f(a) = p, f(b) = q$ in Case 1), and $f(a) = f(b) = k$ in Case 2). 
      
      (iii) For all $i\in \{0, 1, \dots , 2m\}$, $|f^{-1}(i)| = x_i$.
       
       \medskip
       
       Here, we carry the induction on the quantity $m$. For $m\leq 4$ we can see by a direct check that the claim holds. If $m\geq 5$ then there exists $i\in \{0, \dots , m-2\}$ such that none of the triples $(x_{2i}, x_{2i+1}, x_{2i+2}), (x_{2i+2}, x_{2i+3}, x_{2i+4})$ belongs to $B\cup C$. Then using Lemma 1 we reduce the question to the case of a sequence of length $2m-2$.     
       
       \medskip
       
       Now, conditions (i)-(ii) imply that $\Omega _f = S$. But, by interpolation formula [1], for any piecewise linear function one can find a polynomial $p:[0,1]\to \mathbb{R}$ with same range and the same set of local maximums and local minimums. Then $\Omega _p = \Omega _f.$ \ $\square $
      
    \medskip
    
    {\bf Acknowledgement:} I thank Azer Akhmedov for his encouragements and for so many helpful comments and suggestions on the original draft.
    
    \vspace{1cm}
    
    REFERENCES:
    
    \bigskip
    
    [1] \ Whittaker, E. T. and Robinson, G. "Lagrange's Formula of Interpolation." §17 in The Calculus of Observations: A Treatise on Numerical Mathematics, 4th ed. New York: Dover, pp. 28-30, 1967. 
     
    \vspace{1cm}
    
    E-mail: 147414@west-fargo.k12.nd.us
  
  \end{document}